\documentclass[11pt]{article}

\catcode`\@=11
\def\marginnote#1{}

\newcount\hour
\newcount\minute
\newtoks\amorpm
\hour=\time\divide\hour by60 \minute=\time{\multiply\hour by60
\global\advance\minute by-\hour}
\edef\standardtime{{\ifnum\hour<12
\global\amorpm={am}%
        \else\global\amorpm={pm}\advance\hour by-12 \fi
        \ifnum\hour=0 \hour=12 \fi
        \number\hour:\ifnum\minute<10
        0\fi\number\minute\the\amorpm}}
\edef\militarytime{\number\hour:\ifnum\minute<10
0\fi\number\minute}

\def\draftlabel#1{{\@bsphack\if@filesw {\let\thepage\relax
   \xdef\@gtempa{\write\@auxout{\string
      \newlabel{#1}{{\@currentlabel}{\thepage}}}}}\@gtempa
   \if@nobreak \ifvmode\nobreak\fi\fi\fi\@esphack}
        \gdef\@eqnlabel{#1}}
\def\@eqnlabel{}
\def\@vacuum{}
\def\draftmarginnote#1{\marginpar{\raggedright\scriptsize\tt#1}}

\def\draft{\oddsidemargin -.5truein
        \def\@oddfoot{\sl preliminary draft \hfil
        \rm\thepage\hfil\sl\today\quad\militarytime}
        \let\@evenfoot\@oddfoot \overfullrule 3pt
        \let\label=\draftlabel
        \let\marginnote=\draftmarginnote
   \def\@eqnnum{(\theequation)
   \rlap{\kern\marginparsep\tt\@eqnlabel}%
\global\let\@eqnlabel\@vacuum} } 
\relax

\draft
\input{amssym.def}
\input{amssym}
\usepackage{eufrak}
\usepackage{color}

\newtheorem{theorem}{Theorem}

\newtheorem{definition}[theorem]{Definition}

\newtheorem{lemma}[theorem]{Lemma}
\newtheorem{proposition}[theorem]{Proposition}
\newtheorem{remark}[theorem]{Remark}

\def\qq{q^{-1}}
\def\De{\Delta}

\def\gl{\mathfrak{gl}}
\def\A{\cal{A}}
\def\de{\delta}
\def\der{\frac{d}{d\,v}}

\def\gg{\mbox{$\frak g$}}
\def\ot{\otimes}
\def\ott{\odot}
\def\h{{\hbar}}

\def\C{{\Bbb C}}

\def\AA{{\cal A}}
\def\Sym{{\rm Sym\, }}

\def\End{{\rm End}}

\def\vv{V^{\otimes 2}}

\def\lhq{\ifmmode {\cal L}(q,\hbar)\else ${\cal L}(q,\hbar)$\fi}

\def\lqh{\ifmmode {\cal L}(q,\hbar)\else ${\cal L}(q,\hbar)$\fi}

\def\gh{\gg_{\h}}

\def\Tr{{\rm Tr}}

\def\al{{\alpha}}

\def\be{\begin{equation}}
\def\ee{\end{equation}}

\begin{document}

\makeatletter
\renewcommand{\theequation}{{\thesection}.{\arabic{equation}}}
\@addtoreset{equation}{section} \makeatother

\title{Generalizations of Poisson structures related to rational Gaudin model}

\author{
\rule{0pt}{7mm} Dimitri
Gurevich\thanks{gurevich@univ-valenciennes.fr}\\
{\small\it LAMAV, Universit\'e de Valenciennes,
59313 Valenciennes, France}\\
\rule{0pt}{7mm} Vladimir Rubtsov\thanks{volodya@univ-angers.fr}\\
{\small\it UNAM, LAREMA  UMR 6093 du CNRS, Universit\'e d'Angers,
49045 Cedex 01, Angers, France and}\\
{\small\it ITEP, Theory Division, 25, Bol. Tcheremushkinskaya, 117259, Mosocow, Russia}\\
\rule{0pt}{7mm} Pavel Saponov\thanks{Pavel.Saponov@ihep.ru}\\
{\small\it
National Research University Higher School of Economics,}\\
{\small\it 20 Myasnitskaya Ulitsa, Moscow 101000, Russia}\\
\rule{0pt}{7mm} Zoran \v{S}koda\thanks{zskoda@irb.hr}\\
{\small\it Rudjer Bo\v{s}kovi\v{c} Institute,
Zagreb, Croathia}}

\maketitle

\begin{abstract}
The Poisson structure arising in the Hamiltonian approach to the rational Gaudin model looks very similar to  the so-called modified Reflection Equation Algebra.
Motivated by  this analogy, we realize a braiding of the mentioned Poisson structure, i.e. we  introduce a "braided Poisson" algebra associated with an involutive  solution to
the quantum Yang-Baxter equation.  Also, we exhibit  another generalization of the Gaudin type Poisson structure by replacing the first derivative in the current parameter, entering the so-called local form of this structure, by a higher order derivative. Finally,  we introduce  a structure, which combines  both generalizations.
Some commutative  families  in the corresponding  braided Poisson algebra are found.

\end{abstract}

{\bf AMS Mathematics Subject Classification, 2010:} 17B37, 17B80, 81R50, 81R12


\section{Introduction}

This note arises from the following observation. The  Poisson bracket which
enters  the construction of the rational Gaudin model (see \cite{Gau, T}) looks quite similar
to so-called modified Reflection Equation (RE) algebra.   In order to make this parallelism more transparent, we  present this
Poisson bracket in the {\em local} form\footnote{This means that we consider the bracket $\{f(u), g(v)\}$ restricted to the diagonal $u=v$, i.e. for the same value of the current parameter, and similarly for all derivatives of the functions $f$ and $g$  (see (\ref{fam})).}.
In this form  the
Poisson bracket in question (called below {\em Gaudin type Poisson bracket}) consists of a family  of relations between all order
derivatives of the matrix $L(u)$. The  usual flip $P$ is  involved in these relations.
Motivated by the aforementioned analogy, we consider  braided Poisson algebras which are defined by a similar formula but with  a
braiding $R$ instead of the flip $P$.

Let us precise that by a {\em braiding} we mean a solution $R:\vv\to\vv $ to the quantum Yang-Baxter equation
\be
R_{12}R_{23}R_{12}=R_{23}R_{12}R_{23},\,\,{\rm where}\,\, R_{12}=R\ot I,\,  R_{23}=I \ot R,
\label{YB}
\ee
$V$ is a finite--dimensional vector space over the ground field $\C$, and $I$ is the identity operator.

We call $R$ a {\em Hecke symmetry} if it satisfies a complementary condition
$$(R-q\,I)(R+\qq I)=0,\,\, q\in  \C$$
 provided $q\not=1$  and  an {\em involutive symmetry} provided $q=1$.
Below, in order to concretize what type symmetry we are dealing with, and  by slightly abusing the language, we  speak about the Hecke case or the involutive one.

To any braiding $R$ we associate the  Reflection Equation (RE) algebra, generated by entries of a matrix $L=\|l_i^j\|_{1\leq i,j \leq n}$ subject to
the system of relations
\be R\, L_1\, R\, L_1-L_1\, R\, L_1\, R= \h (R\, L_1-L_1\, R). \label{RE} \ee
(We call $L$ the {\em generating} matrix of this algebra.)
More precisely, we use this terminology  for the  algebra (\ref{RE}) provided $\h=0$ and we call it {\em modified} RE algebra
 provided $\h\not=0$.
In terms of the RE  algebra in its two forms (modified and non-modified) we  construct a {\em braided analog}\footnote{We use the term "braided" in the sense of the paper \cite{GS}.
The famous RTT algebra, for instance, is not a braided one according to the definition from this paper.} of
the Poisson structure entering the Gaudin  model and exhibit  a family of elements in involution in the spirit
of the classical theory. The crucial role in our construction is played  by the so-called $R$-trace $\Tr^R L$ whose properties look like those
of the usual (super-)trace. In order to ensure the existence of this trace we assume $R$ to be skew-invertible (see section 4).

Besides,  we introduce another generalization of the Gaudin model by replacing  the first derivative in the current parameter $v$ which enters the corresponding Poisson bracket, realized  in its local form, by a higher  order derivative (section 3).
While this derivative is of order 2 we also exhibit this Poisson structure in a {\em global} form (similar to the standard form of the Gaudin
type Poisson bracket).
We find also a commutative family of elements of the corresponding Poisson algebra generated by the entries of the matrix $L(v)$ and its derivatives on $v$.

In the last section of this note we consider  the structures which combine the both aforementioned  generalizations, namely the braided Poisson brackets with higher order derivatives in its local form. This consideration is preceded  by exhibiting {\em braided Poisson structures} associated with symmetries without parameter (section 4). The aim of this section is to compare braided Poisson structures arising from involutive
symmetries and those arising from the Hecke ones.
In particular, we want to show that the properties of the latter structures are much more  complicated
than those of the former ones. We conclude that  a universal definition
of  braided Poisson structures is somewhat unsustained.
The results of this section enable us to
realize the mentioned braiding of Poisson structures from sections 2 and 3.
We realize
such a braiding by means of a skew-invertible involutive symmetry. We complete
this note with a discussion on an eventual generalization of this construction to the Hecke case.

{\bf Acknowledgements }
This paper has arisen from numerous discussions between participants of Franco-Croatian joint program "COGITO", the project N 24829NH "Syst\`emes Int\'egrables et Structures  Non-commutatives" supported
by EGIDE. Z.\v{S}. and P.S. acknowledge hospitality of LAREMA (Angers) and LAMAV (Valenciennes) respectively, D.G. and  V.R. are grateful to Rudger Bo\u{s}kovi\`{c} Institute (Zagreb) and to prof. S. Meljanac for excellent working conditions and numerous discussions during their visits to Croatia in the framework of the Project. V.R. was  also partially supported by the ANR "DIADEMS", Programme Blan - 01020 and by the RFBR grant 12-01-00525.
The work of one of the authors (P.S.) was carried out within the research grant 13-05-0006 under the National Research University
Higher School of Economics Academic Fund Program support in 2013 and was also partially supported by the RFBR grant
11-01-00980-a. We are grateful to T. Skrypnyk for very useful discussions.

\section{Poisson algebras related to Gaudin model}

In this section we consider Poisson  algebras  related to the rational Gaudin model.  Let
\be L(v)=\|l_i^j(v)\|_{1\leq i, j\leq n} \label{LL} \ee
be  a $n\times n$ matrix
with the entries $l_i^j(v)$ which are meromorphic as function of
 $v\in \C$. Hereafter, $i$ (resp., $j$) is the label
of the line (resp., of the column) to which the element  $l_i^j(v)$ belongs.

Let us define a Poisson bracket\footnote{This bracket differs  from the conventional one by  sign. This change is motivated by our desire to make it more
similar to the defining relations (\ref{RE}) of the modified RE algebra.}
  by
\be \{L_1(u),L_2(v)\}=\left[L_1(u)+ L_2(v), \frac{P}{u-v} \right] \label{odin} \ee
where $L_1=L\ot I, L_2=I\ot L$, and $P$ is the usual flip in $\vv$. Also,  the notation $\{A, B\}$, where $A=\|A_i^j\|$ and
$B=\|B_i^j\|$ are two square matrices of the same size, stands for the matrix with the entries  $\sum_k\{A_i^k, B_k^j\}$.

Taking in consideration that $L_2=P\,L_1\, P$ we can represent the right hand side of (\ref{odin}) as
$$\left[\frac{L_1(u)-L_1(v)}{u-v},\, P\right]. $$

Being written via the  entries of $L(u)$, the relation (\ref{odin}) reads
\be \{l_i^j(u), l_k^l(v)\}= \frac{(l_i^l(u)-l_i^l(v))\de_k^j}{u-v}-\frac{(l_k^j(u)- l_k^j(v))\de_i^l }{u-v}. \label{dva} \ee

The Jacobi identity
\be \{\{L_1(u),L_2(v)\}, L_3(w)\}+\{\{L_3(w),L_1(u)\}, L_2(v)\}+\{\{L_2(v),L_3(w)\}, L_1(u)\}=0 \label{Jac} \ee
is ensured by the fact that the element $r(u)=\frac{P}{u}$ is a classical r-matrix, i.e. it  satisfies the classical Yang-Baxter equation:
\be  \left[r_{12}(u), r_{23}(v)\right]+ \left[r_{12}(u), r_{13}(u+v)\right]+\left[r_{13}(u+v), r_{23}(v)\right]=0. \label{CYB} \ee

If  $f$ (resp., $g$) is a polynomial in generators $l_i^j(u)$ (resp., $l_i^j(v)$), the bracket $\{f,g\}$
can be computed via the Leibniz rule.

Now, we put $u=v+h$  and expand the matrix $L(v+h)$ in the Taylor series  in $h$
\be L(v+h)=\sum_{k=0}^\infty \frac{L^{(k)}(v)\, h^k}{k!} \label{Tay} \ee
where $L^{(k)}(v)=\frac{d^k}{dv^k} L(v) $ stands for the $k$-order derivative of $L(v)$. Comparing the terms containing   the same powers of $h$, we  get
\be \{L^{(k)}_1(v), L_2(v)\}=\frac{\left[L^{(k+1)}_1(v), P\right]}{k+1} \label{case} \ee

By  deriving this equality in $v$ and by using the Leibniz rule for the derivative $\der$
 \be \der\{L^{(k)}_1(v), L_2(v)\}=\{L^{(k+1)}_1(v), L_2(v)\}+\{L^{(k)}_1(v), L^{(1)}_2(v)\}, \label{lei}\ee
we get the following relation
$$\{L^{(k)}_1(v), L^{(1)}_2(v)\}=\frac{\left[L^{(k+2)}_1(v), P\right]}{(k+1)(k+2)}.$$

By continuing this procedure, we recurrently  arrive  to the formula
\be \{L^{(k)}_1(v), L^{(l)}_2(v)\}=\left[L^{(k+l+1)}_1(v), P\right] \al_{1}(k,l), \label{fam} \ee
where
\be \al_{ 1}(k,l)=\frac{k!\, l!}{(k+l+1)!},\,\, k,l=0,1,2,... \label{alph} \ee
(as usual, we assume that $L^{(0)}(v)=L(v)$).

Observe that the coefficients $ \al_{1}(k,l)$ entering  this formula are symmetric:
$\al_1(k,l)=\al_1(l,k)$.   Moreover, the elements
\be \beta_1(k,l,m)=\al_{ 1}(k,l)\al_1(k+l+1,m)=\frac{k!\, l!\, m!}{(k+l+m+2)!} \label{beta} \ee
are invariant with respect to the cyclic permutations of $k,l,m$. This property  ensures  the Jacobi relation for the bracket presented
in the {\em local} form (\ref{fam}).

Thus, we have a family of the relations (\ref{fam}) labeled by couples of naturals  $(k,l)$. This family
with $\al_{1}(k,l)$ given by (\ref{alph}) is equivalent to (\ref{odin}). The passage back to the form (\ref{odin}) can be also
done via the Taylor series (\ref{Tay}).

The  Poisson structure in its local form (\ref{fam}) is defined on the algebra  $\A$ 
of polynomials in the entries of $L(v)$ and all its derivatives. The underlying vector space of this algebra is graded:
its component of degree $k$ is just the span of the elements $\frac{d^k}{d\,v^k}l_i^j(v)$. We extend this grading  by "linearity"
to a graded algebra structure on  $\A$ setting $deg(ab)=deg(a)+deg(b)$. This graded structure on $\A$ is "shifted"
by the Poisson bracket (\ref{fam}) : the components with labels  $k$ and $l$ tmap to that with the label $k+l+1.$

\begin{remark} \label{rem1}
Note that the Poisson bracket
\be \{L_1(v), L_2(v)\}=\left[L_1(v), P\right] \label{case1} \ee
which looks like that (\ref{fam}) with $k=l=0$ but without any derivative in the right hand side,  is nothing but the Lie-Poisson bracket which corresponds to
 the current Lie algebra $\widehat{\gl(n)}$.  Indeed,  (\ref{case1}) is a matrix form  of the following family of  relations
 \be \{l_i^j(v), l_k^l(v)\}= l_i^l(v)\de_k^j-l_k^j(v)\de_i^l . \label{dvaa} \ee
 
However, this bracket does not admit a "natural" extension to the higher derivatives of the type
$\{L_1^{(k)}(v), L_2^{(l)}(v)\}$. 
Indeed, any such extension  \be \{L^{(k)}_1(v), L^{(l)}_2(v)\}=\left[L^{(k+l)}_1(v), P\right]  \label{fam0} \ee
 is not compatible with the action of derivations on the relation (\ref{case1}). Nevertheless, if we disregard this compatibility and treat $k$ and $l$ as labels only, the  Poisson (and the corresponding Lie) algebra structure is well defined.
 \end{remark}
\begin{remark} 
The Lie algebra, defined by formula (\ref{fam}) is bigger than  the current algebra
$\widehat{\gl(n)}$. Indeed,  our algebra is defined on the graded vector space while the underlying space of the  current algebra has the unique component of degree 0. However,  the latter can be converted into the affine algebra, being extended by the
Kac-Moody cocycle.  We have not succeeded to a similar cocycle on our algebra (in this study the local form is very useful). 
It seems that the problem of  constructing an analogue of  quantum affine algebras in the spirit of
\cite{RS} is not well consistent in our setting. 

 \end{remark}

Now  consider a  specialization  of the matrix $L(u)$ of the form
\be L(v)=C+\sum_{p=1}^N A(p) f_p(v) \label{spec} \ee
where $A(p),\,\,p=1,2,...,N$ are matrices with the entries $a_i^j(p)$ subject to (\ref{dvaa}) (where the parameter $v$ is canceled)
for any $p$ and such that
\be \{a_i^j(p), a_k^l(q)\}=0, \,\,\forall i,j,k,l\,\,{\rm if}\,\, p\not=q. \label{comr} \ee
Also, $C$ is a constant matrix. Consequently, its entries Poisson commute with the entries of $A(p)$ for any $p$.

In other words, we have a Poisson bracket $\{\,,\,\}_{\mathfrak{G}}$ where
\be \mathfrak{G}=\mathfrak{gl}(n)\oplus ...\oplus  \mathfrak{gl}(n)\oplus \mathfrak{gl}(n)_{0}, \label{gg}\ee
is  a direct sum.
 Hereafter, by $\{\,,\,\}_{\gg}$ we denote the linear Poisson-Lie bracket corresponding to the Lie algebra $\gg$, and $\gh$ stands for the Lie algebra which
 differs from
 $\gg$ by the factor $\h$ introduced at the Lie bracket. Thus, in the Lie algebra $\mathfrak{gl}(n)_{0}$ 
 the bracket is trivial. 

 Whereas the functional factors $f_p(v)$ are to be found.
More precisely, we want the specification (\ref{spec}) to satisfy   the defining relations (\ref{odin}) or (what is the same) the family of the relations (\ref{fam}).

It is easy to  see that these relations are equivalent to the following differential equation on the factors $f_p(v)$
$$f_p(v)^2=\der f_p(v),\,\, \forall\,p=1,2,...,N.$$
General solution to this equation is $f_p(v)=\frac{1}{v_0-v}$ with any fixed value $v_0$. (It differs by the sign  from the factors in the usual Gaudin model,
see footnote 3.)
Consequently,  the matrix (\ref{spec}) with functional factors $f_p(v)=\frac{1}{v_p-v}$, where $v_p, \, p=1,2,...$, are
any fixed complex numbers (poles), is subject to  the relations (\ref{fam}) provided $\al_{1}$ are defined  by (\ref{alph}). Or, equivalently, the matrix (\ref{spec})
satisfies the relation (\ref{odin}).

The basic property of the Poisson structure above consists in the fact that
\be \{\Tr (L(u))^k, \Tr (L(v))^l\}=0,\,\, \forall\, k,l=0,1,2... \label{inv} \ee
Namely, by using this property for $k=l=2$ and by considering the specialization (\ref{spec})
 one constructs a commutative family of Hamiltonians for the Gaudin model.

More precisely, by considering the quantities
$$\Tr L(v)^2 = \sum_{p=0}^N \frac{\Tr A(p)^2}{(v_p-v)^2} +\sum_p\frac{H(p)}{v_p-v},$$
one gets a family of the quadratic Hamiltonians
$$ H(p) = \Tr C\,A(p) + 2\sum_{j\neq p}\frac{\Tr A(p) A(j)}{v_j - v_p}, 1\leq p \leq N.$$
which commute with each other $[H(i), H(j)] = 0$ for any $1\leq i,j \leq N.$
For detail the reader is referred to  (\cite{CRT, T}).

\section{Generalization of Gaudin type Poisson algebras via  higher order derivatives in $v$}

Below,  we exhibit certain Poisson structures  generalizing those considered in the previous section.
Let us fix an integer $r\geq 2$ and considering  (\ref{fam}) as a pattern, define a Poisson bracket via
\be \{L^{(k)}_1(v), L^{(l)}_2(v)\}=\left[L^{(k+l+r)}_1(v), P\right] \al_{r}(k,l) \label{famq} \ee
where
\be \al_{r}(k,l)=\frac{(k+r-1)!\, (l+r-1)!}{(k+l+2r-1)!(r-1)!},\,\, k,l=0,1,2,... \label{alphq} \ee

The Jacobi identity for such bracket is straightforward to check. This checking is based on the fact that
the terms
 \be \beta_r(k,l,m)=\al_{r}(k,l)\al_p(k+l+r,m)=\frac{(k+r-1)!\, (l+r-1)!\, (m+r-1)!}{(k+l+m+2r-1)!((r-1)!)^2} \label{betaq} \ee
are invariant with respect to the cyclic permutations of $k,l,m$.

Our choice of the coefficients (\ref{alphq}) is motivated by the fact that the matrix $L(v)$ in this case also admits a specialization of the form (\ref{spec})
but with other functional factors $f_p(v)$, namely, those $f_p(v)=\frac{1}{(v_p-v)^r}$. Also, formula (\ref{famq}) with such coefficients is compatible with
its differentiations  in $v$.

Denote  the Poisson  brackets (\ref{famq}) by $\{\,,\,\}_r$. Thus,  the bracket (\ref{fam}) can be treated as a particular case
of  (\ref{famq}) with $r=1$.

Now, we consider the following problem: if it is possible to represent the Poisson bracket $\{\,,\,\}_r$ in a way similar to  the bracket
(\ref{odin})? In other words, we want to  compute the bracket $\{L_1(u), L_2(v)\}$ which is {\em global} (i.e. similar to the standard form  of the bracket (\ref{odin})).  Here we give an answer to this question for $r=2$.

\begin{proposition} The bracket $\{\,,\,\}_2$ can be represented as follows
\be \{L_1(u),L_2(v)\}_2=\left[L_1(u)+ L_1(v), \frac{P}{(u-v)^2}\right]-2 \left[\int_v^u L_1(t)dt,\, \frac{P}{(u-v)^3} \right] \label{odinn} \ee
\end{proposition}

{\bf Proof} can be done in the same way, namely, by setting $u=v+h$ and expanding $L(u)$ in the Taylor series in $h$.
By doing so, we get
$$\sum_{k=0}^\infty\{L_1^{(k)}(v),L_2(v)\}\frac{h^k}{k!}=
\sum_{k=0}^\infty \left[L_1^{(k+2)}(v),P\right]\al(k,0)\frac{h^k}{k!}=$$ $$
\sum_{k=0}^\infty \left[L_1^{(k+2)}(v),P\right]\frac{(k+1)h^k}{(k+3)!}=
\sum_{k=0}^\infty \left[L_1^{(k+2)}(v),P\right]h^k\left(\frac{1}{(k+2)!}-\frac{2}{(k+3)!}\right)=$$ $$
\left[\frac{L_1(v+h)-L_1(v)-L_1^{(1)}(v)h}{h^2},P\right]-2\left[\frac{F_1(v+h)-F_1(v)-F_1^{(1)}(v)h-
F_1^{(2)}(v)h^2\,2^{-1}}{h^3},P\right], $$
where $F(v)$ is a primitive of $L(v)$. It is not difficult to see that the last expression is equivalent
to the right hand side of (\ref{odinn}).

\begin{remark} Nevertheless, formula (\ref{odinn}) differs  from  that (\ref{odin})--no classical r-matrix enters formula (\ref{odinn}).
Thus, a direct verification of the fact that it defines a Poisson bracket  indeed, becomes difficult. However, such a verification for the bracket, realized
under the local form (\ref{famq}), is straightforward.
\end{remark}

Our next aim is to exhibit a family of elements in involution in the Poisson algebra equipped with the bracket  $\{\,,\,\}_{r}$.
 First, we should precise that, similarly to the construction of the previous section, this bracket is well defined on the commutative algebra $\AA$ generated by the
 entries of the matrix $L(v)$ and  its derivatives.

 \begin{proposition}\label{four} In this Poisson algebra the following commutation relations take place
 \be \{\Tr (L(v))^k, \Tr (L(v))^l\}_r=0,\,\, \forall k,l. \label{inv1} \ee
 \end{proposition}

 {\bf{Proof}}\,\,  First, observe that the Poisson brackets $\{\,,\,\}_{r}$ in consideration are of the following form
 \be \{L_1(v), L_2(v)\}_r=[C_1(v), P] \label{matC}\ee
 where $C(v)$ is a matrix, which  up to a factor equals  $L^{(r)}(v)$.
 Consequently, its entries belong to the algebra $\AA$.

Applying  the Leibniz rule to the bracket (\ref{matC}) we get (below we omit the current parameter $v$ and the subscript $r$)
 \be \{ (L_1)^k, (L_2)^l\}
= \sum_{i,j} (L_1)^i\, (L_2)^j\, (C_1\,P - P\, C_1) (L_1)^{k-i-1}\, (L_2)^{m-j-1}.
\label{Poi} \ee
Hereafter in this section we assume that all sums are taken over $0\leq i\leq k-1,\,\,0\leq j\leq l-1$.

Opening the middle brackets in (\ref{Poi}), we transform the first term as follows
 \be  \sum (L_1)^i\, (L_2)^j\, C_1\,P\,  (L_1)^{k-i-1}\, (L_2)^{m-j-1}=\sum (L_1)^i\,  C_1\,P\,(L_1)^j\,  (L_1)^{k-i-1}\, (L_2)^{m-j-1} \label{last} \ee

Here, we have used that $(L_2)^j\, P=P\, (L_1)^j$. Also, we  used the fact that the entries of  $(L_2)^j$ and those of $C_1$ commute with each other. So, we can apply  the first claim  of the following lemma.
\begin{lemma} Let $A$ and $B$ be two square matrices of the same size and the entries  of $A$ commute with those of $B$. Then
$A_1 B_2=B_2 A_1$. Also,  $\Tr\,AB=\Tr\, BA$.  \end{lemma}

For the same reason we can represent the right hand side of  formula (\ref{last}) as
\be \sum (L_1)^i\,  C_1\, (L_1)^{m-j-1}\, P\,(L_1)^j\,  (L_1)^{k-i-1}. \label{od} \ee

 In a similar way we can transform the second term in the right hand side of formula (\ref{Poi})
 $$\sum_{i,j} (L_1)^i\, (L_2)^j\,  P\, C_1\, (L_1)^{k-i-1}\, (L_2)^{m-j-1}=\sum_{i,j} (L_1)^i\,  P\, (L_1)^j\, C_1\, (L_1)^{k-i-1}\, (L_2)^{m-j-1}=$$
 \be \sum_{i,j} (L_1)^i\, (L_1)^{m-j-1} \, P\, (L_1)^j\, C_1\, (L_1)^{k-i-1}. \label{dv} \ee

 Now, apply the operator $\Tr_{12}=\Tr_1\ot \Tr_2$ to the both sides of the relation (\ref{Poi}). On the left hand side we get
 $$ \Tr_{12}\,\{(L_1(v))^k, (L_2(v))^l\}=\{ \Tr_1\,(L_1(v))^k, \Tr_2\,(L_2(v))^l\}=\{ \Tr\,(L(v))^k, \Tr\,(L(v))^l\}.$$
 On the right hand side of (\ref{Poi}) by employing (\ref{od}) and (\ref{dv}) we get
  $$ \Tr_{12} \, \sum \left( (L_1)^i\,  C_1\, (L_1)^{m-j-1}\, P\,(L_1)^j\,  (L_1)^{k-i-1}-(L_1)^i\, (L_1)^{m-j-1} \, P\, (L_1)^j\, C_1\, (L_1)^{k-i-1}\right). $$

 Taking in consideration that $\Tr_2\,P=I$ we transform this formula  to
\be  \Tr \, \sum \left( (L_1)^i\,  C_1\, (L_1)^{m-j-1}\, (L_1)^j\,  (L_1)^{k-i-1}-(L_1)^i\, (L_1)^{m-j-1} \, (L_1)^j\, C_1\, (L_1)^{k-i-1}\right). \label{alm} \ee

Now, by using the second claim of the above lemma we get that the right hand side of (\ref{alm}) vanishes. This implies the claim  of the proposition.

\begin{remark} It seems  that in our  setting  the  stronger identity (\ref{inv})  does not hold.
By contrary,  we are able to show (using the same method as above) that for any natural $m\geq 1$ the following commutation relation holds
 \be \{\Tr (L^{(m)}(v))^k, \Tr (L^{(m)}(v))^l\}=0,\,\, \forall k,l,\, m=1,2,... \label{inv2} \ee

 Whereas, in the frameworks of the initial Gaudin model  by taking the derivative in  $u$ (resp.,  $v$) $m$ (resp., $n$) times of the relation (\ref{inv}) and by setting $u=v$,
 we get a more large family of elements in involution computed at the same value of the current parameter.
\end{remark}

\section{Braided  structures: comparing involutive and Hecke cases}

In this section we compare some structures related to involutive symmetries and Hecke ones. First, remind some facts
about braidings and symmetries (for the detail the reader is referred to \cite{GPS}).

Consider a skew-invertible braiding $R:\vv\to \vv$, where $V$ is a finite dimensional vector space $\dim V=n$.
Let us recall that a braiding $R:\vv\to\vv$ is called {\em skew-invertible} if there exists an operator $\Psi:\vv\to\vv$ such that
$$
\Tr_2 R_{12}\,\Psi_{23}= \Tr_2 \Psi_{12}\,R_{23}=P_{13}.
$$
For any skew-invertible braiding $R:V^{\otimes 2}\to V^{\otimes 2}$ we define two operators $B, C:V \to V$  as follows
\be
B=\Tr_1\Psi,\,\, C=\Tr_2\Psi.
\label{BC}
\ee
We need these operators for introducing braided analogs of  pairings and traces. Recall the corresponding definition.

Let $V^*$ be the dual space to $V$. If $R$ is a skew-invertible braiding, then there exists a unique extension of $R$ up to a braiding
\be
\vv\stackrel{R}{\to} \vv,\quad    V\ot V^*  \stackrel{R}{\to}   V^*\ot V,\quad
 V^*\ot V\stackrel{R}{\to}  V\ot V^*,\quad  (V^*)^{\ot 2}\stackrel{R}{\to}  (V^*)^{\ot 2},
\label{braid}
\ee
such that the pairing
\be
 <\,,\,>: V\ot V^*\to \C,\quad  <x_i,  x^j> = \de_i^j \label{pa}
\ee
is $R$-{\em invariant}. (Here $\{x_i\}_{1\le i\le n}$  and $\{x^j\}_{1\le j\le n}$ are the dual bases in $V$ and $V^*$ respectively\footnote{The basis $\{x^j\}$
satisfying the relation (\ref{pa}) is sometime called {\em right dual}.}.)

The $R$-{\em invariance} of the paring means, by definition, that the
following properties take place
$$
R\langle \,,\,\rangle_{12}=\langle\,,\,\rangle_{23} R_{12} R_{23}\quad {\rm on}\quad  V\ot V^*\ot U,
$$
$$
R\langle \,,\, \rangle_{23}=\langle\,,\,\rangle_{12} R_{23} R_{12}\quad{\rm on}\quad U\ot V\ot V^*
$$
where $U=V$ or $U=V^*$,  and  we assume  $R$ to act on the spaces $U\ot \C$ and  $\C\ot U$ as the usual flip.
In the same sense we speak about the  $R$-invariance of other operators.

Given an operator $B$ (\ref{BC}), we define the pairing $V^*\ot V\to \C$ (where $V^*$ is located on the left hand side of $V$) by the relation
$$
<x^j, x_i>_B=B_i^j,
$$
where $B=\|B_i^j\|$ is the $n\times n$ matrix of the operator $B$ in the basis $\{x_i\}$. The pairing $<\,,\,>_B$ is
also $R$-invariant. Note that the operator $B$ is invertible (see \cite{GPS}).

As for the operator $C$, it is used in the definition of the $R$-trace $\Tr^R A$ where
 $A$ is an arbitrary  $n\times n$ matrix with the entries from the
algebra $\A$. We put
$$
\Tr^R A=\Tr\, (C\cdot A)
$$
where $C = \|C_i^j\|$ is the $n\times n$ matrix of the operator $C$ in the basis $\{x_i\}$ and $\Tr$ is the usual trace.

Let us consider the space $W=V\ot V^*$ and introduce an operator
$$
R_W: W^{\ot 2}\to W^{\ot 2},\quad R_W= R_{23}R_{12}R_{34}R_{23},
$$
where $R$ is the braiding (\ref{braid}). It is easy to see that this operator is also a braiding.

Note that if  $R$ is a Hecke symmetry, the operator $R_W$ is  not a symmetry (either
involutive or Hecke). By contrary, if $R$ is an involutive symmetry, $R_W$ is also so.

Let $L=\|l_i^j\|_{1\leq i,j, \leq n}$  be the generating  matrix of the RE algebra  without the current parameter. The space $W$
can be identified with the linear envelope ${\rm span}(l_i^j)$ of the entries of the matrix $L$ via the following map
\be
W\to {\rm span}(l_i^j): \quad  x_i\ot x^j \mapsto  l_i^j.
\label{ident}
\ee
Thus, the above braiding $R_W$ can be pushed forward to the space $({\rm span}(l_i^j))^{\ot 2}$.

This identification and the pairing $<\,,\,>_B$ above enable us to introduce an $R$-invariant
product in the space ${\rm span}(l_i^j)$. On the basis elements we define it by the rule
$$
l_i^j\circ_B l_k^l\stackrel{{\rm \tiny def}}{=} x_i\otimes\langle x^j,x_k\rangle_B\otimes x^l = x_i\otimes x^lB_k^j =
l_i^l\, B_k^j.
$$
Besides, by introducing the $R$-invariant action $l_i^j: V\to V$
$$
l_i^j(x_k)\stackrel{{\rm \tiny def}}{=}x_i\otimes \langle x^j,x_k\rangle_B = B_k^j x_i,
$$
we can identify the space ${\rm span}(l_i^j)$ with the algebra $\End(V)$.

Now, define a pairing $({\rm span}(l_i^j))^{\ot 2}\to \C$ by composing the above product $\circ_B$ and the
paring (\ref{pa})
\be
<l_i^j,  l_k^l>\stackrel{{\rm \tiny def}}{=} \langle l_i^j\circ_B l_k^l\rangle =
\langle x_i, x^l\rangle B_k^j = \de_i^l\, B_k^j.
\label{l-pair}
\ee
The  paring (\ref{l-pair}) (also  denoted $<\,,\,>_B$) is $R$-invariant too.

Next, discuss the problem of  introducing a  braided analog of the symmetric algebra $\Sym(\gl(n))$. It is tempting to
define a braided analog of the symmetric algebra  of the space $W$ by
\be
T(W)/\langle {\rm Im}(I-R_W)\rangle,
\label{symm}
\ee
where $T(W)$ is the free tensor algebra of the space $W$ and $<I>$ is the ideal generated by a subset $I\subset T(W)$.

However,  this algebra  does not possess {\em a good deformation property}. This means that dimensions of
homogeneous components of this algebra differ from the classical ones  $\dim \Sym^k(\frak{gl}(n))$, provided
 $R$ is a Hecke symmetry which is a deformation of the usual flip (i.e. the Hecke symmetry $R=R(q)$ in question
depends on $q$ and turns into the usual flip $P$ as $q=1$).

In our approach the role of such a "braided symmetric algebra" is played by the non-modified RE algebra (\ref{RE})
($\h=0$). (Here we use the above identification (\ref{ident}).) In contrast with the algebra (\ref{symm}), the RE algebra
has the good deformation property, i.e. for a generic $q$ dimensions of its homogeneous components are classical.
Whereas the modified RE algebra (\ref{RE})  ($\h\not=0$) is treated to be a braided analog of the enveloping algebra $U(\gl(n)_\h)$. It is motivated by the fact that if $R$ is a Hecke symmetry, deforming the flip, its limit as $q\to 1$ is just
the enveloping algebra $U(\gl(n)_\h)$. Besides, it has many other properties similar to  those of the algebra $U(\gl(n)_\h)$ (see \cite{GPS}).

Observe  that if $R$ is a Hecke symmetry, the non-modified RE algebra ($\h=0$) is isomorphic to the corresponding
modified RE algebra ($\h\not=0$). Their isomorphism is established by the following map
\be
L\mapsto \h I-(q-\qq) L,\qquad l_i^j\mapsto \h \de_i^j-(q-\qq)l_i^j, \,\,\,\, q\not=\pm 1.
\label{change}
\ee
Thus, in this case the "braided symmetric" and the "braided enveloping" algebras do not differ from each other.

Nevertheless, the quotients of these two versions of the RE algebras over the ideal generated by the elements
$\Tr^R L$  (which are central in both algebras)
are not isomorphic to each other. One of them turns into the algebra $\Sym(\frak{sl}(n))$, the other one into that $U(\frak{sl(n)}_\h)$
as $q\to 1$.

Also, note that there exists an operator $Q: ({\rm span}(l_i^j))^{\ot 2}\to ({\rm span}(l_i^j))^{\ot 2}$ such that the modified RE algebra (\ref{RE}) can be cast in the following form
\be l_i^j\ot l_k^j-Q(l_i^j\ot l_k^l)=\h[l_i^j, l_k^l] \label{brenv} \ee
where $[\,,\,]:({\rm span}(l_i^j))^{\ot 2}\to {\rm span}(l_i^j)$ stands for a linear operator, which is treated to be a braided analog of the Lie bracket.  The explicit form of this bracket is
\be [\,,\,]=\circ_B (I-Q). \label{lie}\ee

The properties of the corresponding {\em braided Lie algebra} and its enveloping algebra (\ref{brenv}) were studied in \cite{GPS} (similar algebras corresponding to involutive symmetries were introduced  by one of the authors in the 80's, see \cite{G}). We use them
below in order to define the corresponding braided Poisson structures.

Now, let us consider the involutive case. While $R$ is an involutive symmetry, the situation simplifies.

\begin{proposition} If $R$ is an involutive symmetry, the operators $R_W$ and $Q$ become equal to each other. Moreover, the non-modified  RE algebra ($\h=0$)
and the algebra (\ref{symm}) coincide.
 \end{proposition}

Besides, the isomorphism (\ref{change}) between two versions of the RE algebra (modified one and non-modified one) fails.

Moreover, assuming $R$ to be involutive, it is easy to give an axiomatic definition of a braided commutative algebra.

\begin{definition}
An associative algebra $A$ endowed with an involutive symmetry $R:A^{\ot 2} \to A^{\ot 2}$ is called braided commutative if
$$
\circ (a\ot b)= \circ R(a\ot b),\,\,\forall a, b\in A
$$
where $\circ: A^{\ot 2}\to A$ is the product in this algebra which  is assumed be be $R$-invariant. If, besides $A$ is unital,
we also assume that $R(1\ot a)=a\ot 1,\, \forall a\in A$.
\end{definition}

Also, there is a natural definition of a braided Poisson structure on such an algebra.

\begin{definition}
\label{def:brPois}
Let $A$ be a braided commutative algebra in the sense of the previous definition. We say that an $R$-invariant operator
$\{\,,\,\}:  A^{\ot 2}\to A$ is a braided Poisson bracket, if the following axioms are fulfilled

1. $\{a,b \}=-\{\,,\,\}\,R(a\ot b),$

2. $\{a,bc \}=\{ a,b\}c +\{\,,\,\}_{23}R_{12}(a\ot b \ot c),$

3. $\{\,,\,\}\{\,,\,\}_{12}(I+R_{12} R_{23}+R_{23} R_{12})(a\ot b \ot c)=0,\,\,\,\forall\, a,b,c\in A.$
\end{definition}

Let us exhibit two examples of such structures. One of them is the linear bracket (\ref{lie}) extended on the whole
algebra defined by (\ref{brenv}) with $\h=0$ via the braided Leibniz rule, i.e. the property 2 above (we have only to
replace the braided Lie bracket by the corresponding braided Poisson one). The second example arises upon replacing
$l_i^j$ in the right hand side of (\ref{lie}) by $\de_i^j$. Also, note that these brackets are compatible, so any their linear
combination is also an example of such a braided Poisson structure.

These two examples can be generalized  to the Hecke case. However, in this case the axioms 1-3 above are not valid any more. We are only able to write their analogs on the generators and exhibit the way of extending the bracket on whole algebra in question.
 This extension is usually defined via a version of the Leibniz rule.
 Since the Leibniz rule is related to the coproduct in the algebra in question,
we first describe this coproduct in the modified RE algebra. While $\h=1$ this coproduct acts on the generators as
follows\footnote{This coproduct was deduced in \cite{GPS} from  the braided bi-algebra structure
discovered in the (non-modified) RE algebra by Majid (see \cite{M}).}
\be
\De(1)=1\ot 1,\,\,\, \De(l_i^j)=l_i^j\ot 1+1\ot l_i^j-(q-\qq)\sum_k l_i^k\ot l_k^j
\label{cop}
\ee
On the whole modified RE algebra the coproduct (\ref{cop}) must be extended with the use of the braiding $R_W$.
For instance,
$$\De(l_i^j\, l_k^l)=\De(l_i^j)\De(l_k^l)=((l_i^j)_1\ot (l_i^j)_2)((l_k^l)_1\ot (l_k^l)_2)=(l_i^j)_1\, \widetilde{(l_k^l)_1}\ot \widetilde{(l_i^j)_2}(l_k^l)_2$$
where  $\widetilde{(l_k^l)_1}\ot \widetilde{(l_i^j)_2}=R_W((l_i^j)_2\ot (l_k^l)_1)$ (here we use the Sweedler's notation).

If $R$ is an involutive symmetry (i.e. $q=1$), this coproduct becomes similar to the classical one. On the generators it
takes the usual form $\De(l_i^j)=l_i^j\ot 1+ 1\ot l_i^j$. Its extension to the whole algebra in question must be done  via the braiding $R_W=Q$ (see proposition 7).

Now, we introduce a convenient matrix notation, which enables us to cast all considered structures in a form useful
for braiding of the Gaudin type Poisson brackets.

Let $R:\vv\to\vv$, $\dim V=n$ be a Hecke or involutive symmetry and $L=\|l_i^j\|_{1\leq i,j,\leq n}$ be a matrix. We put
\be
L_{\overline{1}}=L_1,\quad
L_{\overline{2}}=R_{12}\, L_1\, R_{12}^{-1},\quad
\dots,\quad L_{\overline{i+1}}=R_{i\,i+1}\, L_{\overline{i}}\,R_{i\,i+1}^{-1}.
\label{def}
\ee
If $R$ is involutive, we can replace $R^{-1}$ by $R$.

Also, below we use the  notation $A\ott B$ for the product of two squared  matrices of the same size but with the
entries  multiplied in the sense of the tensor product:
$$
(A\ott B)_i^j=\sum_k A_i^k \ot B_k^j.
$$

In this matrix notation the defining relations of the RE algebra can be written as
$$
R\,(L_{\overline{1}}\ott L_{\overline{2}})-(L_{\overline{1}}\ott L_{\overline{2}}\, ) R=\h(L_{\overline{2}}-L_{\overline{1}}),
$$
or, equivalently, as
\be
L_{\overline{1}}\ott L_{\overline{2}}-R^{-1}(L_{\overline{1}}\ott L_{\overline{2}}\, ) R=\h(L_{\overline{1}}R^{-1}-R^{-1}L_{\overline{1}}).
\label{new2}
\ee
If $R$ is involutive this relation is also equivalent to
\be L_{\overline{1}}\ott L_{\overline{2}}-L_{\overline{2}}\ott
L_{\overline{1}}=\h(L_{\overline{1}}- L_{\overline{2}})R.
\label{new1}
\ee

If $R$ is a Hecke symmetry, the action of the operator $R_W$ can be written as follows
$$
R_W(L_{\overline{1}}\ott L_{\overline{2}})=L_{\overline{2}}\ott L_{\overline{1}},
$$
whereas for the operator $Q$ we get
$$
Q(L_{\overline{1}}\ott L_{\overline{2}})=R^{-1}(L_{\overline{1}}\ott L_{\overline{2}}) R.
$$

Following the pattern arising from formula (\ref{new2}) (where we put $\h=1$), we define the linear Poisson bracket
on the generators by the rule
\be
\{L_{\overline{1}}, L_{\overline{2}}\}=L_{\overline{1}}R^{-1}-R^{-1}L_{\overline{1}}.
\label{defi}
\ee
(For the notation $\{A, B\}$ see section 2.)

The properties of this bracket are similar to those of the braided Lie algebras considered in \cite{GPS}. However,
one needs to complement formula (\ref{defi})  with a rule for extending
the bracket to monomials of higher order.
In our current setting  the Leibniz rule does not have a universal form similar to that from definition 9. However, we define the extension of the bracket
via the coproduct (\ref{cop}).

Thus, we have
\be
\{L_{\overline{1}},L_{\overline{2}}\ott L_{\overline{3}} \}=\{L_{\overline{1}},L_{\overline{2}}\}\ott L_{\overline{3}}
+L_{\overline{2}}\ott\{L_{\overline{1}}, L_{\overline{3}} \}-(q-\qq)\{L_{\overline{1}},L_{\overline{2}}\}\ott
\{L_{\overline{1}}, L_{\overline{3}} \}.
\label{leib}
\ee
Note that $\{L_{\overline{1}},L_{\overline{3}}\}= \{L_{\overline{1}},R_{23}L_{\overline{2}}R_{23}^{-1}\}=R_{23}
\{L_{\overline{1}},L_{\overline{2}}\}R_{23}^{-1}$.

In a similar manner the brackets
$$\{L_{\overline{1}}\ott L_{\overline{2}}\}, L_{\overline{3}}, \,\,\, \{L_{\overline{1}}\ott L_{\overline{2}}, L_{\overline{3}}\ott L_{\overline{4}}\}$$
and so on can be computed.

Whereas, the braided analog of the Jacobi relation is similar to that from \cite{GPS}
$$\{L_{\overline{1}}, \{L_{\overline{2}}, L_{\overline{3}}\}\}=\{\{L_{\overline{1}}, L_{\overline{2}}\}, L_{\overline{3}}\}+
\{L_{\overline{2}}, \{L_{\overline{1}}, L_{\overline{3}}\}\}.$$

 The last property we want to mention, is the skew-symmetry: in the Hecke case it takes the form
 $$ \{\,,\,\}((q^2+q^{-2})(L_{\overline{1}}\ott L_{\overline{2}})+R^{-1}\, (L_{\overline{1}}\ott L_{\overline{2}}) R+
 R\, (L_{\overline{1}}\ott L_{\overline{2}}) R^{-1})=0.$$
 Note that the expression in the bracket is the image of the symmetrization operator on $W^{\ot 2}$ (see \cite{GPS})  realized in the matrix form.

To conclude, we note that to define a Poisson structure on a braided algebra is a somewhat subtle deal. Especially,
the form of  the "braided Leibniz rule" is not a priori clear (see the discussion at the end of the next section). However,
in the involutive case this rule (expressed by the property 2 above) is simple enough and via the usual $R$-matrix
technique (see \cite{GPS} and the references therein) it entails the formula
$$
\{(L_{\overline{1}})^k, (L_{\overline{2}})^l\}=\sum_{i=0}^{k-1}\sum_{j=0}^{l-1} (L_{\overline{1}})^i\ott
(L_{\overline{2}})^j\ott  \{L_{\overline{1}}, L_{\overline{2}}\}\ott (L_{\overline{1}})^{k-i-1}\ott (L_{\overline{2}})^{l-j-1},
$$
which is used in the next section.

\section{Braided Poisson structures of the Gaudin type}

First, consider the braided version of the Gaudin type bracket (\ref{odin}). To this end we replace the usual flip $P$
in that formula by a skew-invertible involutive symmetry $R$. Then we get
\be
\{L_{\overline{1}}(u), L_{\overline{2}}(v)\}=\left[L_{\overline{1}}(u)+ L_{\overline{2}}(v), \frac{R}{u-v} \right].
\label{odinnn}
\ee
Hereafter, the notation $L_{\overline{i}}(u)$ stands for the matrices defined according to formula  (\ref{def})  where $R$ is
the given involutive symmetry  (without the current parameter). Thus, the braiding and differentiations in the parameter do not affect each
other.

It is not difficult to see that the relation (\ref{CYB}) is fulfilled with $r(u)=\frac{R}{u}$. Here the fact that $R$ is involutive is
crucial. However, we prefer to deal with the corresponding braided Poisson  structure by casting it in the  local form similar
to (\ref{fam}).

Let us consider the algebra $\A$ generated by the entries of the matrix $L(v)$ and its derivatives which are subject to the relations
$$
L^{(k)}_{\overline{1}}(v)\ott  L^{(l)}_{\overline{2}}(v)=L^{(l)}_{\overline{2}}(v)\ott L^{(k)}_{\overline{1}}(v),\,\,\forall k,\, l.
$$
This algebra looks like that from definition 8 but now it becomes graded. We also treat it as a {\em braided commutative} algebra.   We assume that the bracket (\ref{odinnn}) is defined on the algebra

Our next aim is to define  the bracket (\ref{odinnn}) in the local form on this algebra.

On expanding $L_{\overline{1}}(u)=L_{\overline{1}}(v+h)$ in the Taylor series in $h$ (in analogy with section 2), we
conclude that formula (\ref{odinnn}) is equivalent to the family
\be
\{L^{(k)}_{\overline{1}}(v), L^{(l)}_{\overline{2}}(v)\}=[L^{(k+l+1)}_{\overline{1}}(v),\, R]\, \al_1(k,l)
\label{famm}
\ee
with coefficients $ \al_1(k,l)$ given by formula (\ref{alph}).

In a similar manner we can realize a braiding of the Poisson structures considered in section 3. Namely, we set
\be
\{L^{(k)}_{\overline{1}}(v), L^{(l)}_{\overline{2}}(v)\}_r=\left[L^{(k+l+r)}_{\overline{1}}(v), R\right] \al_{r}(k,l)
\label{famqq}
\ee
with the coefficients $\al_{r}(k,l)$ defined by formula (\ref{alphq}).

This Poisson structure can also be cast in the form similar to (\ref{odinn}). For $r=2$  we have
\be
\{L_{\overline{1}}(u),L_{\overline{2}}(v)\}_2=\left[L_{\overline{1}}(u)+ L_{\overline{1}}(v), \frac{R}{(u-v)^2}\right]-
2 \left[\int_v^u L_{\overline{2}}(t)dt,\, \frac{R}{(u-v)^3} \right].
\label{odi}
\ee

The properties of all these braided brackets are similar to those exhibited definition 9 but now the algebra is graded and the bracket is compatible with the gradation.

Consider the braided  enveloping algebra similar to (\ref{new1}), which is generated by the entries of the matrices $L^{(k)}(v),\, k=0,1,2,...$ subject to the  relations
\be
R\, L_1^{(k)} (v)\, R\, L_1^{(l)}(v) - L_1^{(l)} (v)\, R\, L_1^{(k)}(v)\, R=\al_r(k,l) ( R\, L_1^{(k+l+r)}(v)-
L_1^{(k+l+r)}(v)\, R)
\label{REA}
\ee
where  the coefficients $\al_r(k,l)$ are defined by formula (\ref{alph}).
(The corresponding Lie bracket can be also  readily defined.)

It is not difficult to see that this algebra becomes a braided bi-algebra being equipped with the coproduct
\be \De L^{(k)}(v)=L^{(k)}(v)\ot 1+ 1\ot L^{(k)}(v)\label{copro}\ee

 extended to the algebra $\A$ by means of the operator
$$
Q(L^{(k)}_{\overline{1}}(v)\ott  L^{(l)}_{\overline{2}}(v))=L^{(l)}_{\overline{2}}(v)\ott L^{(k)}_{\overline{1}}(v).
$$
Namely, this property leads to a proper Leibniz rule in the braided Poisson algebra in question.

Now, we  discuss a specialization of the matrix $L(v)$ similar to (\ref{spec}). Let $\frak{gl}(R)$ be the braided
Lie algebra defined in the space ${\rm span}(a_i^j)$ by (\ref{defi}) but with the Lie bracket instead of the Poisson
one and the matrix $A=(a_i^j)$ instead of $L$.

By following the pattern (\ref{gg}) we consider a braided analog of the direct sum (\ref{gg}) but  without the last component
$$
\mathfrak{G}(R)=\mathfrak{gl}(R)\oplus_R ...\oplus_R  \mathfrak{gl}(R).
$$
This means that the commutation relations in the algebra $\mathfrak{G}(R)$ are
\be
[A_{\overline{1}}(p),A_{\overline{2}}(q)]:=A_{\overline{1}}(p)\ott A_{\overline{2}}(q)-A_{\overline{2}}(p)\ott A_{\overline{1}}(q)=(A_{\overline{1}}R-R A_{\overline{1}})\de(p,q), \label{abc}
\ee
where  $\de(p,q)$ is the Kronecker symbol and $p$ and $q$ are labels of the components (recall that $R$ is involutive).

\begin{proposition}
The matrix $L(v)=\sum_{p=1}^N A(p) f_p(v)$, where the entries of the matrices $A(p)$ belong to the braided Lie algebra $\mathfrak{G}(R)$ (i.e. they are subject to (\ref{abc})) fulfills the relations (\ref{odinnn}) iff the functional factors are $f_p(v)=\frac{1}{v_p-v}$.
Also, the matrix $L(v)=\sum_{p=1}^N A(p) f_p(v)$ is subject to  (\ref{famqq}) iff the functional factors are $f_p(v)=\frac{1}{(v_p-v)^r}$.
\end{proposition}

In the same manner as the classical result (see \cite{T}) the following proposition can be proven.

\begin{proposition}
The Hamiltonians
$H(p)=\sum_{j\not=p}\frac{\Tr^R A(p) A(j)}{ v_j-v_p}, \, 1\leq p\leq N$, where the matrices $A(p)$  are subject to the relations
(\ref{abc})  commute with each other in the sense of the braided Poison bracket $\{\,,\,\}_{\mathfrak{G}(R)}$.
\end{proposition}

Also, in the spirit of the proposition \ref{four} we can prove the following.

\begin{proposition} If the matrices $L^{(k)}(v)$ are subject to the bracket (\ref{famqq}) then the quantities
$\Tr^R  (L(v))^{k}$, $k=0,1,2,...$ commute with each other  in the sense of the bracket (\ref{famqq}) extended to the whole algebra via the coproduct (\ref{copro}).
\end{proposition}

We want only to emphasize the crucial points of the proof. First, formula (\ref{Poi}) is still valid, provided $L_i(v)$ is replaced by $L_{{\overline{i}}}(v)$. Since the matrix $C(v)$ equals up to a factor a derivative of the matrix $L(v)$ the relation
$$
L_{{\overline{1}}}(v)\ott C_{{\overline{2}}}(v)=C_{{\overline{2}}}(v)\ott L_{{\overline{1}}}(v)
$$
is valid.

Now, in order to get a braided analog of the proposition 4, it suffices to apply the following lemma.

\begin{lemma}
If two $n\times n$ matrices $A$ and $B$  are subject to the relation
\be
A_{{\overline{1}}}(v)\ott B_{{\overline{2}}}(v)=B_{{\overline{2}}}(v)\ott A_{{\overline{1}}}(v)
\label{rela}
\ee
then
\be
\Tr^R AB=\Tr^R BA.
\label{tr-comm}
\ee
\end{lemma}

{\bf Proof} We will show that this claim is valid even in the Hecke case. Let us rewrite the relation (\ref{rela}) as
\be
A_1\, R\, B_1= R\, B_1\, R\, A_1\, R^{-1}
\label{trr}
\ee
and apply the $R$-trace $\Tr^R_{12}=\Tr^R_{1}\ot \Tr^R_{2}$ to the both sides of this equality.
Next, we apply the relation
$$
\Tr^R_{12}\,  R\, X_{12} \, R^{-1}=\Tr^R_{12} X_{12},
$$
where $X_{12}\in \End(\vv)\ot \A$ is an endomorphism of the space $\vv$ with coefficients belonging to  any associative  algebra $\A$.
If we set $X_{12}=B_1\, R\, A_1$ and use the property $\Tr^R_{2}R_{12} = I_1$ we come to the result (\ref{tr-comm}).
\medskip

Completing the note, we want to discuss the following problem: is it possible to define an analogous structure related to
a Hecke symmetry. Let $R$ be such a symmetry and fix an integer $r\ge 1$. Consider an associative algebra generated by the entries
of the matrices $L^{(k)}(v)$, $k=0,1,2,...$ subject to the  relations (\ref{REA}) with the chosen Hecke symmetry $R$. It is not difficult
to define a braided Lie bracket (in the spirit of \cite{GPS}) so that the algebra (\ref{REA}) becomes the enveloping algebra of the
corresponding braided Lie algebra. However, it is not clear whether there exists a coproduct endowing it with a braided bi-algebra
structure. This is reason why we do not know any consistent form of the "braided Leibniz rule" in the corresponding "braided Poisson
structure". Also, the deformation property of this algebra is not clear. Nevertheless, in this associative algebra the problem of finding an
analog of the Cayley-Hamilton identity in the spirit of \cite{T} is of interest.

\end{document}